\documentclass[11pt,letterpaper]{amsart}
\usepackage[utf8]{inputenc}
\usepackage{verbatim}
\usepackage{amsmath,amssymb,amsfonts,amsthm}
\usepackage{bigints}
\usepackage{hyperref}
\usepackage{array}
\numberwithin{equation}{section}

\usepackage{verbatim}
\usepackage{amsmath,amssymb,amsfonts,amsthm}
\usepackage{bigints}
\usepackage{hyperref}
\usepackage{array}

\usepackage{amsthm}

\usepackage{mathtools}

\mathtoolsset{showonlyrefs}

\newtheorem{definition}{Definition}[section]

\newtheorem{theorem}[definition]{Theorem}

\newtheorem{lemma}[definition]{Lemma}

\newtheorem{remark}[definition]{Remark}

\newcommand{\norm}[1]{ \left\lVert#1\right\rVert}

\newcommand{\spec}{\operatorname{spec}}

\newcommand{\lipstandard}{\operatorname{Lip}}

\thanks{
This work is supported by
AP23487589, 
 PID2023-150984NB-I00 funded by MICIU/AEI/10.13039/501100011033/ FEDER, EU, 
the CERCA Programme of the Generalitat de Catalunya and the Severo Ochoa, and Mar\'ia de Maeztu
Program for Centers and Units of Excellence in R\&d (CEX2020-001084-M).
M. Saucedo is supported by  the Spanish Ministry of Universities through the FPU contract FPU21/04230.
 S. Tikhonov is supported
by 2021 SGR 00087. }

\author{Miquel Saucedo}
\address{M.  Saucedo,  Centre de Recerca Matemàtica\\
Campus de Bellaterra, Edifici C
08193 Bellaterra (Barcelona), Spain}
\email{miquelsaucedo98@gmail.com }

\author{Sergey Tikhonov}
\address{S. Tikhonov, Centre de Recerca Matem\`{a}tica\\
Campus de Bellaterra, Edifici C
08193 Bellaterra (Barcelona), Spain;
ICREA, Pg. Lluís Companys 23, 08010 Barcelona, Spain,
 and Universitat Autònoma de Barcelona.}
\email{ stikhonov@crm.cat}

\subjclass[2010]{Primary  
42B05, 42A20; Secondary 46E35.}
\keywords{Fourier coefficients,
Absolute convergence, Smooth function spaces}

\title[Kahane-Katznelson-de Leeuw  theorem]{Kahane-Katznelson-de Leeuw  theorem   %in smooth spaces
and absolute convergence of Fourier series
}

\begin{document}
\begin{abstract}
   We extend the Kahane-Katznelson-de Leeuw theorem 
   to smoothness spaces by showing that 
for any $g \in W^{l,2}(\mathbb{T}^d)$, there exists a function $f\in C^l(\mathbb{T}^d)$ satisfying $|\widehat{f}(n)|\geq |\widehat{g}(n)|$ and $$\omega_r(D^l f,t)_\infty \approx \omega_r(D^l g,t)_2, \quad t>0. $$
 We apply this result %As application, %  of this result,   we
 to solve the Bernstein problem of finding necessary and sufficient conditions for the absolute convergence of multiple Fourier series.
 Finally, we explore the absolute integrability  of Fourier transforms.
  %We also discuss related results for Fourier transforms.
 
\end{abstract}

\maketitle
% \tableofcontents
\section{Introduction and main results}
\subsection{Kahane-Katznelson-de Leeuw theorem for moduli of smothness}
The classical Kahane-Katznelson-de Leeuw theorem \cite{kkdl} states that for any $g\in L^2(\mathbb{T}^d)$ there exists a function $f\in C(\mathbb{T}^d)$ with $\norm{f}_\infty \approx \norm{g}_2$ such that $$|\widehat{f}(n)|\geq |\widehat{g}(n)|.$$ 
There are several generalizations  \cite{Kis24,nazarov} of this result to spaces different from the space of continuous functions. 
 %For instance, among many other results, 
In particular, Kislyakov \cite{kislyakov88,Kis24} proved the following variant of this theorem for the $C^s$ spaces.
Recall that 
$$\|f\|_{C^s(\mathbb{T}^d)}=
\|f\|_{C(\mathbb{T}^d)}+
\max_{|\alpha|=s}\Big\|\partial^\alpha f
\Big\|_{C(\mathbb{T}^d)}
$$
and $\mathbb{N}_0=\mathbb{N}\cup \{0\}.$
%$f\in C^s(\mathbb{T}^d)$ if $\partial_\alpha f\in C(\mathbb{T}^d)$,  $|\alpha|=l$.

\begin{theorem}[\protect{\cite[Th. B]{Kis24}}]
\label{theorem:KKdL}
    Let $(c_k)_{k\in \mathbb{Z}^d}$ satisfy
    \begin{equation}
        \sum_{k \in \mathbb{Z}^d} \big(c_k |k|^s\big)^2 \leq 1.
    \end{equation} Then, for any $s\in \mathbb{N}_0$ there exists a function $f$ such that
    \begin{itemize}
    \item[$(i)\,\,\,\,$]\quad $\widehat{f}(n)\in \mathbb{R};$
               \item[$(ii)\,\,\,$]\quad $\norm{f}_{C^s(\mathbb{T}^d)}\lesssim 1;$
        \item[$(iii)\,$]\quad $|\widehat{f}(n)|\geq |c_n|, \, n \in \mathbb{Z}^d$.   
    \end{itemize}
\end{theorem}
%\begin{proof}
 We note that even though the condition $\widehat{f}(n)\in \mathbb{R}$ is not stated, it is straightforward to see from the techniques in \cite{kislyakov88} that it can be assumed.
%\end{proof}

Our first main result, Theorem \ref{th:main}, is an % says that the 
 analogue of Theorem \ref{theorem:KKdL} for moduli of smoothness. 
 
\begin{definition} Let $l\in \mathbb{N}_0$,
$r\in \mathbb{N}$,
and $q\geq 1$.
\begin{itemize}
    \item[$(i)$] We define the $r$-th modulus of smoothness of the $l$-th derivative of a function $f\in W^{l,q}(\mathbb{T}^d)$ if $ q<\infty$ or $f\in C^{l}(\mathbb{T}^d)$ if $q=\infty$
    %$ \to \mathbb{C}$ 
     by
    \begin{equation}\label{definition-mod}
\omega_r(D^l f,t)_q:=\sup _{\substack{0<|h|<t\\|\alpha|=l
}} {\norm{\Delta^r_h (\partial_\alpha f)}}_q,
    \end{equation}
    where 
    $\Delta_h^1 g(x)=g(x)-g(x+h)$ and $\Delta_h^r=\Delta_h^1(\Delta_h^{r-1}).$
  \item[$(ii)$]  We say that $\omega:[0,1]\to \mathbb{R}_+$ is $r$-quasiconcave if
    \begin{itemize}
        \item[$\cdot$] $\omega(t)$ and $t^r/\omega(t)$ are non-decreasing for $t\in [0,1];$
        \item[$\cdot$] $\omega(0)=0$.
    \end{itemize}
 \item[$(iii)$]
  For % Let $l,r\in \mathbb{N}$ and 
  an $r$-quasiconcave $\omega$, we say that $f\in 
   {\operatorname{Lip}^{r,l}_q(\omega;\mathbb{T}^d)}
   $ if $$\sup _{0<t<1} \frac{\omega_r(D^l f,t)_q}{\omega(t)}<\infty.$$
     \end{itemize}
\end{definition}
We note that for every $f$ there exists an $r$-quasiconcave $\omega$ such that $\omega_r(D^l f,t)_q\approx \omega(t)$, see for instance \cite[Remark 4]{koltikh}. Besides, if there exists  $\varepsilon>0$ such that $\omega(t)/t^\varepsilon$ is almost increasing, i.e.,
$\omega(t_1)/t_1^\varepsilon\lesssim
\omega(t_2)/t_2^\varepsilon$, $t_1\le t_2,$ then 
   $\operatorname{Lip}^{r,l}_q(\omega;\mathbb{T}^d)=  {\operatorname{Lip}^{r+l,0}_q(t^l \omega;\mathbb{T}^d)}$,
   see 
   \cite[Property 11]{koltikh}.

Our first main result is the
Kahane--Katznelson--de Leeuw theorem for moduli of smoothness. % reads as follows:
\begin{theorem}
\label{th:main}
 Let $l\in \mathbb{N}_0$
 and
$r\in \mathbb{N}$.
For any $g \in W^{l,2}(\mathbb{T}^d)$ there exists $f\in C^l(\mathbb{T}^d)$ such that $$|\widehat{f}(n)|\geq |\widehat{g}(n)|$$ and $$\omega_r(D^l f,t)_\infty \approx \omega_r(D^l g,t)_2, \quad t>0. $$
\end{theorem}
Equivalently, this result can be written as follows: 
{\it 
Let $r,l\in \mathbb{N}$ and $\omega$ be $r$-quasiconcave. Then, for any $g\in 
{\operatorname{Lip}^{r,l}_2(\omega;\mathbb{T}^d)}$
there exists a continuous function  $f$ such that}
$|\widehat{f}(n)|\geq |\widehat{g}(n)|$ and 
$$\norm{f}_{\operatorname{Lip}^{r,l}_\infty(\omega;\mathbb{T}^d)}\approx \norm{g}_{\operatorname{Lip}^{r,l}_2(\omega;\mathbb{T}^d)}. $$

\begin{remark}
    From Theorem \ref{th:main} we conclude that 
 $$%   \begin{align*}
        \norm{\widehat{f}(n)\,}_X\lesssim F\big(\omega_r(D^l f,\cdot)_2\big) \mbox{ holds for any  } f\in W^{l,2}(\mathbb{T}^d)
        $$
        if and only if
%        \\
$$\norm{\widehat{f}(n)\,}_X\lesssim F\big(\omega_r(D^l f,\cdot)_\infty\big) \mbox{ holds for any } f\in C^l(\mathbb{T}^d),
$$
   % \end{align*}
    provided that $\norm{\cdot}_X$  and $F(\cdot)$ are order-preserving functionals on sequences and functions, respectively, that is,
    \begin{align}
    \begin{split}\label{lattice}
    |a_n|\lesssim |b_n| \,\,\mbox{ for all }\,\, n\in \mathbb{Z}^d  &\implies \norm{a}_X\lesssim\norm{b}_X;\\
    |a(t)|\lesssim |b(t)|\,\, \mbox{ for all }\,\, t>0 &\implies F({a})\lesssim F({b}).
    \end{split}
        \end{align}
 
\end{remark}

\subsection{Absolute convergence of multiple Fourier series}We apply Theorem \ref{th:main} to study 
the absolute convergence of Fourier series.
We use the notation
$\mathcal{A}_p(\mathbb{T}^d)$ to denote the Wiener space of 
$p$-absolutely
convergent Fourier series, that is, 
$$\mathcal{A}_p(\mathbb{T}^d)=\left\{f\in L^1(\mathbb{T}^d):\sum_{n\in \mathbb{T}^d} |\widehat{f}(n)|^p<\infty \right\}.$$
The problem of determining the smooth spaces contained in $\mathcal{A}_p(\mathbb{T}^d)$ has been investigated since 1914, when Bernstein  proved that\footnote{
In fact, more is true: there is a function from $\lipstandard  1/2$, which
cannot be brought into $\mathcal{A}(\mathbb{T})$ by any correction on a set not of full measure, see \cite{olevskii}.
} % following result:
$$\lipstandard \alpha  \subset \mathcal{A}(\mathbb{T})=\mathcal{A}_1(\mathbb{T})\quad \mbox{ if and only if} \quad\alpha > 1/2.$$

In 1934, Bernstein posed the question of determining when  the Lipschitz space with a ge\-ne\-ral majorant $\omega$  belongs to $\mathcal{A}_1(\mathbb{T})$. In the one-dimensional case the answer to this question (see  \cite{bary, kahane}) is given as follows:
% answer i stated the prob
%This result was later extended to 
%Lipschitz spaces with a general majorant $\omega$, see  \cite{bary, kahane}:
\begin{equation}
\label{condi:bari}
     \operatorname{Lip(\omega;\mathbb{T})}= \operatorname{Lip}^{1,0}_\infty(\omega;\mathbb{T})
\subset \mathcal{A}(\mathbb{T})\;\; \mbox{ if and only if} \;\;\int_0^1\frac{\omega(t)}{t^\frac12} \frac{dt}{t}<\infty.
\end{equation}

In the multivariate case, an investigation  of 
 the Bernstein type result  is much more delicate.
% attracted a lot of attention.
 % for   the 
% Lipschitz spaces is
% not known in the full generality.
We start with the sufficient condition (see \cite{peetre}) 
 given by:

%We now recall the classical Bernstein theorem, which  gives sufficient conditions for the absolute convergence of the Fourier series of a function $f$ function in terms of its smoothness.

%The result 
%As a consequence we are able to 
%obtain necessary and sufficient conditions in terms of the moduli of smoothness
%for the absolute convergence of the Fourier series. This problem has been studied for a long time and there are many available results.

%A first classical result is the Bernstein Theorem, which gives sufficient conditions for the absolute convergence of the Fourier series of a f function in terms of its smoothness.

\begin{theorem} 
\label{thm:previous} 
Define, for $\theta,p>0$ and $q\ge 1,$
        $$B_{q,p}^\theta(\mathbb{T}^d)=\left\{f\in L^q(\mathbb{T}^d): |f|_{B_{q,p}^\theta}=\Big(\int_0 ^1 \Big(t^{-\theta} \omega_{\lceil\theta\rceil+1}(f,t)_q\Big)^p \frac{dt}{t}
        \Big)^\frac1p <\infty\right\}.$$        
       Then,
        $$\displaystyle B_{2,1}^{\frac d2}(\mathbb{T}^d)\subset \mathcal{A}(\mathbb{T}^d).$$

\end{theorem}
In particular, this implies that $ {\operatorname{Lip}(t^{\frac{d}2+ \varepsilon} ;\mathbb{T}^d)}\subset \mathcal{A}(\mathbb{T}^d)$, $\varepsilon>0$, see  \cite[Chapter VII]{stein}.
In 1965, Wainger showed that the smoothness parameter $\frac{d}{2}$  is  sharp, i.e.,
 $ {\operatorname{Lip}(t^{\frac{d}2} ;\mathbb{T}^d)}\not\subset \mathcal{A}(\mathbb{T}^d)$. Equivalently, we have 
%even if $L^2$ smoothness is replaced by $L^\infty$ smoothness. 
%, more precisely,
\begin{theorem}[\protect\cite{wainger65}]
\label{th:wainger}
 Let $l\in \mathbb{N}_0$. Then
        $$%\begin{align}
            C^{l}(\mathbb{T}^{2l})\not \subset \mathcal{A} 
            $$
            and
            $$
             {\operatorname{Lip}^{1,l}_\infty(t^{\frac 12};\mathbb{T}^{2l+1})}
            \not \subset \mathcal{A}.
        $$
\end{theorem}

There have been several attempts (see, e.g., \cite{KasMel23,Kis24, koregrin, szasz, Musielak, nowak,
 timan}) to find sharp conditions for the absolute convergence of Fourier series of functions on $\mathbb{T}^d$
 from  the  Lipschitz space ${\operatorname{Lip}^{r,l}_p (\omega)}$ when $p=2$ or $p=\infty$; however, a complete solution to this problem remains unknown.
  %In more detail,
 
 Here we mention two results. In the odd dimensional case,
 %the following analogue of 
\eqref{condi:bari} was extended 
 to
\begin{equation}
\label{condi:kore}
     {\operatorname{Lip}^{1,l}_\infty(\omega;\mathbb{T}^{2l+1})}
\subset \mathcal{A}(\mathbb{T}^{2l+1})\quad \mbox{ if and only if} \quad\int_0^1\frac{\omega(t)}{t^\frac12} \frac{dt}{t}<\infty,
\end{equation} 
 provided that $\omega(t)/t^{\varepsilon}$ and $t^{1-\varepsilon}/\omega(t)$
 are non-decreasing with 
some  $\varepsilon>0$, see, e.g., \cite{koregrin}.
 
In the even dimensional case,  the authors of \cite{Kas, KasMel23, Kis24} 
have recently investigated the limiting cases that lie between Bernstein's and Wainger's results and  obtained the following

%$\lip_\infty(\mathbb{T}^d)$.
%A complete solution of the problem on absolute convergence of the Fourier series on $\mathbb{T}^d$
%is not known for  $\lip_\infty(\mathbb{T}^d)$.

\begin{theorem}
[\protect\cite{Kis24}]
\label{thm:previous2}
 Let $l\in \mathbb{N}$.
Then
$$\operatorname{Lip}^{1,l}_\infty\big(\big(\log \frac2t\big)^{-\frac{1}{2}};\mathbb{T}^{2l}\big)\not \subset \mathcal{A}.$$   
\end{theorem}

We also mention the following (slightly weaker) result given in 
  \cite{KasMel23}:  $$\displaystyle \operatorname{Lip}^{1,1}_\infty\big(\big(\log \frac2t\big)^{\eta-\frac{1}{2}};\mathbb{T}^2\big)\not \subset \mathcal{A},\quad \eta>0.$$

\iffalse
characterize the $(r,l,d,\omega)$ for which $B_\infty^{r,l\omega}$ is contained in the Wiener algebra $\mathcal{A}$, the space of functions with absolutely convergent Fourier coefficients.

\fi
We now give the characterization of the absolute convergence of the Fourier series in terms of the moduli of smoothness, thus fully answering the question. 

\begin{theorem}
\label{coro:hardy}
 Let 
 $l\in \mathbb{N}_0$,
$r\in \mathbb{N}$,
 and $\omega$ be an $r$-quasiconcave majorant. Let $0<p< 2$. Then, the following are equivalent:
\begin{enumerate}
    \item $
    \operatorname{Lip}^{r,l}_\infty(\omega; \mathbb{T}^{d})
    \subset \mathcal{A}_p;$
      \item $\operatorname{Lip}^{r,l}_2(\omega; \mathbb{T}^{d})
    \subset \mathcal{A}_p;$
    \item one of the following holds:
    \begin{itemize}
       \item[$(i)$] $d(\frac 1p- \frac12)<l,$
       \item[$(ii)$] $l<d(\frac 1p- \frac12)<l+r$ and $$\int_0 ^1 \omega^p(t) t^{-(d(\frac1p -\frac12)-l)p} \frac{dt}{t}<\infty,$$
        \item[$(iii)$] $d(\frac 1p- \frac12)=l$ and $$\int_0 ^1 \frac{\omega^p(t)}{ \big(\log \frac2t\big)^{\frac{p}{2}}} \frac{dt}{t}<\infty.$$
    \end{itemize}
\end{enumerate}
   \end{theorem}
\begin{remark}
\textnormal{(1)} 
Specializing  condition
 (3)(ii) in Theorem
\ref{coro:hardy} for $L^q$--moduli of smoothness, we note that if $l<d(\frac 1p- \frac12)<l+r$, then  for any $q\in [1,\infty]$
$$
\int_0 ^1
 t^{p(-d(\frac1p -\frac12)+l)}\omega_r^p(D^l f,t)_q
 \frac{dt}{t}
 \approx
 \int_0 ^1
 t^{-dp(\frac1p -\frac12)} \omega_{r+l}^p(f,t)_q
\frac{dt}{t}
\approx|f|^p_{B^{d(\frac1p -\frac12)}_{q,p}}.
 %\approx|f|_{B^{d(\frac1p -\frac12)}_{2,p}},
 $$
% and 
% $\|f\|_2 +Q\approx \|f\|_2+|f|_{B^{d(\frac1p -\frac12)}_{2,p}}.$
% cf. condition (3)(ii) in Theorem \ref{coro:hardy}.
Regarding condition (3)(iii),
 we similarly observe that
%, \eqref{maj-log}
$$
\int_0 ^1 \frac{
\omega_r^p(D^l f,t)_q
}{ \big(\log \frac2t\big)^{\frac{p}{2}}} \frac{dt}{t}\approx %<\infty \iff 
\int_0 ^1 \frac{
\omega_1^p(D^l f,t)_q
}{\big(\log \frac2t\big)^{\frac{p}{2}}} \frac{dt}{t} =: %\iff
|f|^p_{ \mathfrak{B}^{l, -\frac 12}_{q,p}},
$$
where $\mathfrak{B}^{l, -\frac 12}_{q,p}$
is the Besov space with logarithmic smoothness.
Various characterizations and properties of both Besov spaces
$B^{l}_{q,p}$ and 
$\mathfrak{B}^{l, m}_{q,p}$ can be found in 
\cite{oscar}. (We emphasize that the logarithmic Besov spaces are defined differently in \cite{oscar}.) In particular, 
$|f|_{ \mathfrak{B}^{l, m}_{q,p}}$ is equivalent to the corresponding functional with $\sup _{|\alpha|=l
} {\norm{\Delta^r_h (\partial_\alpha f)}}_q
$
in place of
$\omega_r(D^l f,t)_q$.

 \textnormal{(2)}  
 Let $0<p<2\leq q \leq \infty$. We say that $X$ is a $(r,l,q)$-lattice if it satisfies $$f\in X, \,\omega_r(D^lg,t)_q\lesssim \omega_r(D^l f,t)_q \implies g \in X.$$ 
 Note  the the condition $\omega_r(D^lg,t)_q\lesssim \omega_r(D^l f,t)_q$
is a contraction condition introduced by Beurling \cite{Beurling}.

 From Theorem \ref{coro:hardy} we deduce that for $2\leq q\leq \infty$ the largest $(r,l,q)$-lattice $X$ contained in $\mathcal{A}_p$ is $$X= \begin{cases}
     B^{d(\frac1p -\frac12)}_{q,p}, &l<d(\frac1p- \frac12)<l+r,\\
     \mathfrak{B}_{q,p}^{d(\frac1p- \frac12), -\frac 12}, &l=d(\frac1p- \frac12).
 \end{cases}
 $$ 
Note that the spaces $\mathfrak{B}_{q,p}^{d(\frac1p- \frac12), -\frac 12}$
and $B^{d(\frac1p -\frac12)}_{q,p}$ are not comparable for $q\ne2$. This follows from \cite[Chapter 12]{oscar}. 

For the case $q=2$, we would like to stress  that 
even though $$  \mathfrak{B}_{2,p}^{d(\frac1p- \frac12), -\frac 12}\subset B^{d(\frac1p -\frac12)}_{2,p} \subset \mathcal{A}_p,$$ the space 
$\mathfrak{B}_{2,p}^{d(\frac1p- \frac12), -\frac 12}
 $
 is still the largest lattice,
 since $B^{d(\frac1p -\frac12)}_{2,p}$ is not a $(r,l,q)$-lattice if $l=d(\frac1p -\frac12).$

 %We also note that even though $  B_{q,p}^{d(\frac1p- \frac12), \frac p2}\subset B^{d(\frac1p -\frac12)}_{q,p} \subset \mathcal{A}_p$, this is not a contradiction because $B^{d(\frac1p -\frac12)}_{q,p}$ is not a $(r,l,q)$-lattice if $l=d(\frac1p -\frac12).$

We note that the Fourier inequalities
$
\|{\widehat{f}(n)}\|_X\lesssim 
\norm{f}_Y$
between rearran\-gement-inva\-ri\-ant lattices $X$ and $Y$ have been recently characterized in
\cite{saucedo2025fouriertransformextremizerclass}.

 \iffalse
 Note that in the case 
 $d(\frac 1p- \frac12)=l$ %, cf. condition (3)(iii),
 we only have
$$|f|^p_{B^{d(\frac1p -\frac12)}_{2,p}}
\lesssim
\int_0 ^1 \frac{
\omega_r^p(D^l f,t)_2
}{ \log_+^{\frac{p}{2}}(t^{-1})} \frac{dt}{t}.
$$
\fi

 \textnormal{(3)}  
The results in Theorems \ref{thm:previous}--\ref{thm:previous2} as well as  \eqref{condi:bari} and \eqref{condi:kore} are particular cases of Theorem \ref{coro:hardy} with $p=1$
for specific choices of parameters. Indeed, Theorem \ref{thm:previous} follows by setting $r=\lceil \frac d2\rceil +1,l=0$; conditions \eqref{condi:bari} and \eqref{condi:kore}, by setting $r=1$ and $d=2l+1$; Theorem \ref{thm:previous2} follows by setting $r=1$, $\omega(t)= \big(\log \frac2t\big)^{-\frac{1}{2}}$ and $d=2l$.
 %and noting that $\int_0^1 \frac{1}{\log_+(t^{-1})} \frac{dt}{t}=\infty$.
Theorem \ref{th:wainger} follows from the other results.
\end{remark}

Throughout the paper, we  use the notation $F \lesssim G$  to mean that $F \le C G$ with a constant $C=C(p,q,r,l,d)$ that may change from
line to line; $F \approx G$ means that both $F \lesssim G$ and $G \lesssim F$ hold. 
Besides, $ B_\infty(r)=\{x:\norm{x}_\infty\leq r\}$ and, as usual, $$\partial_\alpha f=\frac{\partial^\alpha f}{\partial^{\alpha_1}x_1\cdots
\partial^{\alpha_d}x_d}.
$$

\iffalse
\begin{theorem}[Bernstein]
    Assume that $f$ satisfies
    $\int_0 ^1 t^{l-\frac d2}{\omega(t)}\frac{dt}{t}<\infty$
\end{theorem}
Bernstein,

Wainger
\begin{theorem}[Wainger]
    For every $l\in \mathbb{N}$ there exist functions in $C^{l}(\mathbb{T}^{2l})$ and  $C^{l,\frac 12}(\mathbb{T}^{2l+1})$ with non-absolutely convergent Fourier series.
\end{theorem}

Kislyakov Kashin

\fi

\medskip

\section{Auxiliary results}
%\tableofcontents

The next two lemmas are variants of the results from \cite{gogatishvili2003discretization} adapted  to the discrete case.
%The first result in this direction can be found in \cite{oskolkov}.
Note that the construction of the sequence $\mu_k$ in Lemma 
\ref{lemma:discre} appeared already 
 in Oskolkov's paper \cite{oskolkov}.

\begin{lemma}[Discretizing sequences, \protect{cf. \cite[Definition 2.4]{gogatishvili2003discretization}}]
\label{lemma:discre}
   Fix $\lambda>4$ and $r>0$. Let $(\omega_n)_{n=1}^\infty $be such that $\omega_n \searrow 0$ and $\omega_n n^r \nearrow$. Then there exists an increasing sequence $(\mu_k)_{k=1}^{L+1}$ with $L \in [1, \infty]$ such that
    \begin{enumerate}
        \item $\mu_1=1$;
        \item if $L<\infty$, then $\mu_{L+1}=\infty$;
        \item  $\mu^r_{k+1}\geq \lambda \mu^r_k,\,k \in [1,L]$;
        \item $\lambda\omega_{\mu_{k+1}}\leq \omega_{\mu_{k}},\,k \in [1,L]$;
        \item  $\lambda\omega_{\mu_{k}} \mu^r_{k}\leq \omega_{\mu_{k+1}} \mu^r_{k+1},\,k \in [1,L-1]$; 
        \item there are two sets of integers $I$ and $J$ such that  $\mathbb{N}\cap [1,L]=I \cup J$ and 
    \begin{align*}
        &\omega_{\mu_{k}} \approx \omega_{\mu_{k+1}-1}, \qquad k\in I
        \\&\omega_{\mu_{k}} \mu^r_{k} \approx  (\mu_{k+1}-1)^r\omega_{\mu_{k+1}-1}, \qquad  k\in J.
        \end{align*}
        \end{enumerate}
\end{lemma}
\begin{proof}
Set $\mu_1=1$ and define the sequence $\mu$ recursively as follows:
    $$\mu_{k+1}:=\min \Big \{ n: n^r \omega_n > \lambda  \mu^r_k  \omega_{\mu_k}\text{ and } \lambda\omega_n<  {\omega_{\mu_k}}\Big\}.$$ If for $k_0$ the set is empty, then we set $L=k_0$ and we stop.

    The proof of properties (1)--(5) is straightforward. To see that (6) holds observe that if $k\leq L-1$, then for $m=\mu_{k+1}-1$ we either have $m^r\omega_m\leq \lambda   \mu^r_k \omega_{\mu_k} $ or $\lambda\omega_m \geq   \omega_{\mu_k}$. If $k=L$, then any $n \geq \mu_{k}$ satisfies $\mu^r_k\omega_{\mu_k} \leq n^r\omega_n \leq \lambda \mu^r_k \omega_{\mu_k} $.
\end{proof}
The sequence $(\mu_k)_{k=1}^{L+1}$ from
Lemma
\ref{lemma:discre}
will be called {\it discretizing} sequence for $(\omega_n)_{n=1}^\infty $.

% $h_l=\sum_{n=1}^\infty \frac{\alpha_n}{n^r+l^r}$. 
 
 \begin{lemma}[\protect 
 {cf. \cite[Lemma 3.6]{gogatishvili2003discretization}}]
 \label{theorem:split}
Assume that $p,r>0$.
Let $\alpha$ be a non-negative sequence and set
 $$\overline\omega^p_l:=\sum_{n=1}^\infty \frac{\alpha_n}{n^{pr}+l^{pr}},\qquad l\in \mathbb{N}.$$  Let also   $0<q\leq 1$. Then, for any non-negative $f_j$,
\begin{equation}
    \sum_{n=1}^\infty \alpha_n \left( \sum_{j=1}^\infty \frac{f_j}{j^{pr/q}+n^{pr/q}}\right)^ q \approx \sum_{k\in I\cup J} \left( \sum_{\mu_{k} \leq l<\mu_{k+1}} \overline \omega_l ^{\frac{p}{q}} f_l \right)^q,
\end{equation}
where
$\mu_{k}$ is a discretizing sequence for 
$\overline\omega_l^p$.

\end{lemma}
    \begin{proof}
 First, we observe that  $\overline\omega_l^p \searrow 0$ and $l^{pr} \overline\omega_l^p \nearrow $, so by 
 Lemma
\ref{lemma:discre}
 there exists a discretizing sequence $\mu$.
 
        Let us show the $"\lesssim"$ estimate.
        Observe that
         \begin{align}
             &\sum_{n=1}^\infty \alpha_n \left( \sum_{j=1}^\infty \frac{f_j}{j^{pr/q}+n^{pr/q}}\right)^q \leq \left(\sum_{k\in I}+ \sum_{k \in J}\right) \sum_{n=1}^\infty \alpha_n  \left( \sum_{j=\mu_k}^{\mu_{k+1-1}} \frac{f_j}{j^{pr/q}+n^{pr/q}}\right)^q\\
&\lesssim \sum_{k\in I} \sum_{n=1}^\infty \frac{\alpha_n}{\mu^{pr}_k+n^{pr}}   \left( \sum_{j=\mu_k}^{\mu_{k+1-1}} {f_j}\right)^q \\
&\qquad\qquad+ \sum_{k\in J}(\mu_{k+1}-1)^{pr} \sum_{n=1}^\infty \frac{\alpha_n}{n^{pr}+(\mu_{k+1}-1)^{pr}} \left( \sum_{j=\mu_k}^{\mu_{k+1-1}} {f_jj^{-pr/q}}\right)^q
\\         
&\lesssim \sum_{k\in I} \overline\omega^p_{\mu_{k}}  \left( \sum_{j=\mu_k}^{\mu_{k+1-1}} {f_j}\right)^q + \sum_{k\in J} (\mu_{k+1}-1)^{pr}\overline \omega^p_{\mu_{k+1}-1}  \left( \sum_{j=\mu_k}^{\mu_{k+1-1}} {f_jj^{-pr/q}}\right)^q \\
&\approx  \sum_{k\in I\cup J} \left( \sum_{\mu_{k} \leq l<\mu_{k+1}} \overline\omega_l ^{\frac{p}{q}} f_l \right)^q,
         \end{align}
         where the last equivalence follows from item (6) of the previous lemma.

         To obtain the $"\gtrsim"$ part, we first prove that for $k\leq L$

         \begin{equation}
         \label{eq:piece}
             \sum_{j=\mu_{k-1}}^{\mu_{k+2}-1} \frac{\alpha_j}{j^{pr}+n^{pr}}\approx \overline\omega^p_n, \quad \mu_k \leq n <\mu_{k+1},
         \end{equation} (if $k-1\leq1$ or $k+2\geq L+1$, then $j$ starts at $1$ or $j$ ends at $\infty$, respectively). 
         
         For this purpose we estimate
         \begin{align}\label{first}
             \sum_{j=1}^{\mu_{k-1}-1} \frac{\alpha_{j}}{j^{pr}+n^{pr}} & \leq 2 \left(\frac{\mu_{k-1}-1}{n}\right)^{pr}\sum_{j=1}^{\mu_{k-1}-1} \frac{\alpha_{j}}{j^{pr} + (\mu_{k-1}-1)^{pr}} \nonumber\\&\leq  2 \left(\frac{\mu_{k-1}-1}{n}\right)^{pr}\overline\omega^p_{\mu_{k-1}-1} \leq \frac{2}{\lambda}\overline\omega^p_{n}
             \end{align}
             and
\begin{align}\label{second}\sum_{j=\mu_{k+2}}^{\infty} \frac{\alpha_{j}}{j^{pr}+n^{pr}} \leq 2 \sum_{j=\mu_{k+2}}^\infty \frac{\alpha_{j}}{j^{pr} + \mu^{pr}_{k+2}} \leq  2\overline \omega^p_{\mu_{k+2}} \leq \frac{2}{\lambda}\overline\omega^p_{n}, \end{align}
where the last estimates in \eqref{first} and \eqref{second} follow from items (5) and (4) in Lemma \ref{lemma:discre}, respectively.
Hence,
 \begin{equation}
             \sum_{j=\mu_{k-1}}^{\mu_{k+2}-1} \frac{\alpha_j}{j^{pr}+n^{pr}}\geq  (1- \frac 4 \lambda)\overline\omega^p_n, \quad \mu_k \leq n <\mu_{k+1},
         \end{equation} 
and
\eqref{eq:piece} is shown.

         In light of Minkowski's inequality and \eqref{eq:piece}, we have
         \begin{align}
             \sum_{n=1}^\infty \alpha_n \left( \sum_{j=1}^\infty \frac{f_j}{j^{pr/q}+n^{pr/q}}\right)^q &\approx \sum_{k \in I \cup J} \sum_{n=\mu_{k-1}}^{\mu_{k+2}-1} \alpha_n \left( \sum_{j=1}^\infty \frac{f_j}{j^{pr/q}+n^{pr/q}}\right)^q \\
             &\gtrsim   \sum_{k \in I \cup J} \left( \sum_{j=1}^\infty f_j  \left(\sum_{n=\mu_{k-1}}^{\mu_{k+2}-1}  \frac{\alpha_n}{j^{pr}+n^{pr}}\right)^\frac 1q\right)^q\\
             &\geq  \sum_{k \in I \cup J} \left( \sum_{j=\mu_k}^{\mu_{k+1}-1} f_j  \left(\sum_{n=\mu_{k-1}}^{\mu_{k+2}-1}  \frac{\alpha_n}{j^{pr}+n^{pr}}\right)^\frac 1q\right)^q \\
             &\approx \sum_{k \in I \cup J} \left( \sum_{j=\mu_k}^{\mu_{k+1}-1} f_j \overline\omega_j^{\frac pq}\right)^q.
         \end{align}
    \end{proof}
    We now show that every $r$-quasiconcave $\omega$ can be approximated by an $\overline \omega$ satisfying the hypothesis of Lemma \ref{theorem:split}.
\begin{lemma}
\label{lemma:majorant}
    Let $\omega$ be $r$-quasiconcave and $p>0$. Then there exists an $r$-quasiconcave $\overline \omega$ such that
    \begin{equation}
        \omega^p(1/N)\approx \overline \omega^p (1/N) =\sum_{n=1}^\infty \frac{\alpha_n}{N^{pr}+n^{pr}},
    \end{equation} for some non-negative $(\alpha_n)$.
\end{lemma}
\begin{proof}
    First, observe that the function $f(t)=\omega^p(t^{\frac{1}{pr}})$ is non-decreasing and $f(t)/t$ is non-increasing. Thus, the concave majorant of $f$, $$g(x)=\inf_{\substack{f\leq h \\
    \mbox{\tiny concave}}} h(x)$$ satisfies $f\approx g$. Since $g$ is non-decreasing and concave we can write, for some $\beta\geq 0$,
    $$f(t)\approx g(t)  \approx \int_0 ^\infty \min(t,s)\beta(s) ds, \quad 0\leq t \leq 1.$$
    Second, by the definition of $f$ and by a change of variables, for some $\alpha,\gamma\geq 0$ we have
    $$\omega^p(t^{-1}) \approx \int_0^\infty \min(t^{-pr},s^{-pr})\gamma(s) ds \approx \sum_{n=1}^\infty \frac{\alpha_n}{n^{pr}+ t^{pr}},\quad t\geq 1.$$ We complete the proof  by setting
    $\,\displaystyle\overline \omega^p_l= \sum_{n=1}^\infty \frac{\alpha_n}{n^{pr}+ l^{pr}}$.
\end{proof}

    The next result is a Littlewood-Paley type characterization of the Lipschitz spaces using discretizing sequences.
\begin{lemma}

\label{lemma:l2}
    Let $l\in \mathbb{N}_0$ and
$r\in \mathbb{N}$. For an $r$-quasiconcave
$\omega$, we define 
    $\omega_N:=\omega(1/N)$. Let $$f(x)=\sum_{n\in \mathbb{Z}^d} a_n e^{2 \pi i \langle n, x \rangle}\in L^2(\mathbb{T}^d)$$ and
    $$R_n:=\Big(\sum_{\norm{k}_\infty =n} |a_k|^2\Big)^\frac12.$$
    The following are equivalent:
    \begin{itemize}
        \item[$(i)$] we have
        \begin{equation}
    \label{eq:L2func} 
\norm{f}_{\operatorname{Lip}^{r,l}_2(\omega;\mathbb{T}^{d})}
    \lesssim 1;
    \end{equation} 
   \item[$(ii)$] for any $N\in \mathbb{N}$, we have \begin{equation} \label{eq:L2coeff}
        \sum_{n=1}^\infty \Big|R_n n^l \min\Big(1,\frac{n^r}{N^r}\Big)\Big|^2 \lesssim \omega_N^2; %=:\omega^2_N ,
    \end{equation} 
    %where  $R_n^2=\sum_{\norm{k}_\infty =n} |a_k|^2$;
    %and $\omega_N:=\omega(1/N)$;
  \item[$(iii)$]  for all $k \in [1,L], $ we have
    \begin{align}
        \label{eq:L2disc}
        \sum_{n=\mu_k}^{\mu_{k+1}-1} |n^lR_n|^2 \lesssim& \omega_{\mu_k}^2,\quad \mbox{ if }\quad  k\in I \end{align}
      {and}
         \begin{align}  \label{eq:L2disc2}
        \sum_{n=\mu_k}^{\mu_{k+1}-1} |n^{l+r}R_n|^2 \lesssim& ((\mu_{k+1}-1)^r\omega_{\mu_{k+1}-1} )^2, \quad \mbox{ if }\quad k\in J,
    \end{align}
where 
 $(\mu_k)_{k=1}^{L+1}$ is the  discretizing sequence for $(\omega_n)_{n=1}^\infty$ and  $L$ are given in Lemma \ref{lemma:discre}.
\end{itemize}

\end{lemma}
\begin{proof}
    The equivalence \eqref{eq:L2func} $\iff$ \eqref{eq:L2coeff} is well known and follows from the relation
    $$\widehat{\Delta_h^r f }(n)= \left(1-e^{2 \pi i \langle n,h\rangle}\right)^r \widehat{f}(n).$$ 
       Let us now prove that \eqref{eq:L2coeff} implies \eqref{eq:L2disc} and \eqref{eq:L2disc2}.
        Set $$\alpha_n :=R_n n^l.$$ If \eqref{eq:L2coeff} holds,
    setting $N=\mu_k$ we deduce 
     $\sum_{n=\mu_k}^{\mu_{k+1}-1}|\alpha_n^2|  \lesssim \omega^2_{\mu_k}$, which is \eqref{eq:L2disc}. To see \eqref{eq:L2disc2}, we set $N=\mu_{k+1}-1$ in \eqref{eq:L2coeff}.

     We now show that \eqref{eq:L2disc} and \eqref{eq:L2disc2} imply \eqref{eq:L2coeff}.
     Before starting the proof, we observe that \eqref{eq:L2disc} and \eqref{eq:L2disc2} imply
        \begin{align}
    \label{eq:33}
        &\sum_{n=\mu_k}^{\mu_{k+1}-1} |\alpha_n|^2 \lesssim \omega_{\mu_k}^2, \quad k \in I \cup J,\\
        \label{eq:44}
        &\sum_{n=\mu_k}^{\mu_{k+1}-1} |n^r\alpha_n|^2 \lesssim (\mu_{k+1}-1)^{2r} \omega_{\mu_{k+1}-1}^2,\quad k \in I \cup J.
    \end{align}Indeed, if $k\in I$, \eqref{eq:L2disc} coincides with \eqref{eq:33}; if $k\in J$, then \eqref{eq:L2disc2} implies
      $$\sum_{n=\mu_k}^{\mu_{k+1}-1} |\alpha_n|^2 \leq   \sum_{n=\mu_k}^{\mu_{k+1}-1} \big|\frac{n^r}{\mu^r_{k}}\alpha_n\big|^2 \lesssim \frac{1}{\mu_k^{2r}}(\mu_{k+1}-1)^{2r}(\omega_{\mu_{k+1}-1} )^2 \approx \omega_{\mu_k}^2,$$ that is, \eqref{eq:33} holds.
      
Similarly, for $k\in J$, \eqref{eq:L2disc2} is \eqref{eq:44}; if $k\in I$, then \eqref{eq:L2disc} implies
$$\sum_{n=\mu_k}^{\mu_{k+1}-1} |n^r\alpha_n|^2 \leq (\mu_{k+1}-1)^{2r} \omega_{\mu_k}^2 \approx (\mu_{k+1}-1)^{2r} \omega_{\mu_{k+1}-1}^2, $$ that is, \eqref{eq:44} is valid.

All that remains is to show that \eqref{eq:33} and \eqref{eq:44} imply \eqref{eq:L2coeff}. For any $N\in\mathbb{N},$
let $\mu_k \leq N<\mu_{k+1}$ and write
      \begin{align} 
       \left(\sum_{j<k} + \sum_{j=k} +\sum_{j>k}\right) \left(\sum_{n={\mu_j}}^{\mu_{j+1}-1}
    \Big|\alpha_n \min\Big(1, \frac{n^r}{N^r}\Big)\Big|^2
       \right)=:I_1+I_2+I_3.
    \end{align}
To estimate $I_1$, using \eqref{eq:44}, the monotonicity of $\omega_n n^r$ and the property $\lambda \omega_{\mu_k} \mu^r_k \leq \omega_{\mu_{k+1}} \mu^r_{k+1} $ for $k\leq L-1$, we deduce that
    \begin{align}
        I_1&\leq \sum_{j<k} \left(\sum_{n={\mu_j}}^{\mu_{j+1}-1}  \Big|\alpha_n  \frac{n^r}{N^r}\Big|^2\right)\\
        & \lesssim \frac{1}{N^{2r}}\sum_{j<k} (\mu_{j+1}-1)^{2r} \omega^2_{\mu_{j+1}-1} \\
        &\lesssim \sum_{j<k} \frac{\mu^{2r}_{j+1}}{N^{2r}} \omega^2_{\mu_{j+1}}\lesssim  \frac{\mu^{2r}_{k}}{N^{2r}} \omega^2_{\mu_{k}} \leq \omega_N^2.
    \end{align}
    Similarly, for $I_3$, using now \eqref{eq:44},
    \begin{align}
      I_3\leq \sum_{j>k} \left(\sum_{n={\mu_j}}^{\mu_{j+1}-1} |\alpha_n |^2\right)  \lesssim \sum_{j>k}  \omega_{\mu_{j}}^2\lesssim \omega_{\mu_{k+1}} ^2  \leq \omega_N^2.
    \end{align}

    To estimate $I_2$, there are two possibilities:
    First, if $k\in I$, then using \eqref{eq:33} we have
    \begin{align}
        \sum_{n={\mu_k}}^{\mu_{k+1}-1} \Big|\alpha_n \min\Big(1, \frac{n^r}{N^r}\Big)\Big|^2 \leq  \sum_{n={\mu_k}}^{\mu_{k+1}-1} |\alpha_n|^2 \lesssim \omega^2_{\mu_k} \approx \omega^2_N,
    \end{align}
    where the last estimate follows from item (6) of Lemma \ref{lemma:discre}.
%    where the last estimate follows from $k\in I$.
    Second, if $k\in J$, then using \eqref{eq:44} and Lemma \ref{lemma:discre}, we have
    \begin{align}
        \sum_{n={\mu_k}}^{\mu_{k+1}-1} \Big|\alpha_n \min\Big(1, \frac{n^r}{N^r}\Big)\Big|^2&\leq  \sum_{n={\mu_k}}^{\mu_{k+1}-1} |\alpha_n \frac{n^r}{N^r}|^2\\
        &\lesssim \frac{1}{N^{2r}} \Big(\omega_{\mu_{k+1}-1} (\mu_{k+1}-1)^r\Big)^2 \\
        &\approx \frac{1}{N^{2r}} (\omega_{N} N^r)^2=\omega_N^2.
    \end{align}
    %where the last estimate follows from item (6) of Lemma \ref{lemma:discre}.
    The proof is now complete.
\end{proof}
The last lemma relates the norm of a polynomial with those of its derivatives.
\begin{lemma}[Direct and reverse Bernstein's inequalities]
\label{lemma:bernstein}
Let $N,l,r\in \mathbb{N}_0$ with $r\geq 1$. Then,
    \begin{enumerate}
        \item the trigonometric polynomial $
        T_N(x)=\sum_{\norm{k}_\infty \leq N} c_k e^{2 \pi i \langle k , x\rangle}$ satisfies
        \begin{equation*}
            \norm{T_N}_{C^{l+r}} \lesssim N^r \norm{T_N}_{C^l};
        \end{equation*}
        \item the function 
        $F_N(x)=\sum_{\norm{k}_\infty \geq N} c_k e^{2 \pi i \langle k , x\rangle}$ satisfies
        \begin{equation*}
            \norm{F_N}_{C^{l+r}} \gtrsim N^r \norm{F_N}_{C^l}.
        \end{equation*}
    \end{enumerate}
\end{lemma}
\begin{proof}
    The Bernstein inequality for $d=1$ is well known, see  \cite[Th. 3.13]{zygmund_fefferman_2003}. The multidimensional result follows by applying it to each variable separately.

    The
    reverse Bernstein inequality
 is also classical in the one-dimensional case. We give here the proof for the multidimensional case.

    First, by Jackson's inequality (see, e.g., \cite[Property 12]{koltikh}), for any $|\alpha|=l$ and $f$ we have  
\begin{equation}
    \inf \norm{T-\partial_\alpha f}_\infty \lesssim \omega_r(\partial_\alpha f,N^{{-1}})_\infty \lesssim N^{-r} \norm{f}_{C^{l+r}},
\end{equation}
where the infimum is taken over all trigonometric polynomials 
$T$ with the property that 
$\spec(T)\subset B_\infty(N/2)$.

Second, we observe that for any $g$
     \begin{equation}
         \inf_{\spec(T)\subset B_\infty(N/2)} 
         \norm{T-g}_\infty \approx \norm{V^d_{{{N/2}}}(g)-g}_\infty,
     \end{equation}
where $V^d_{{{N/2}}}$ is given by
 %Thus, since we know that, defining 
 \begin{equation}
    \label{eq:VP}
V^d_n(h)=\sum_{\norm{k}_\infty \leq 2n} \widehat{h}(k) e^{ikx} \prod_{i=1}^d\min\left(1, 2-\frac{|k_i|}{n}\right).
    \end{equation}
    Finally, since  $V^d_{N/2}(\partial _\alpha F_N)=0$, we conclude that
     \begin{equation}
         \norm{\partial_\alpha F_N}_\infty \lesssim N^{-r} \norm{F_N}_{C^l+r},
     \end{equation} whence the result follows.
\end{proof}
\medskip
\section{Proofs of Theorems
\ref{th:main} and
\ref{coro:hardy}
}

  \begin{proof}[Proof of Theorem \ref{th:main}]
  % Before starting the proof, note that without loss of generality we may assume that $a_n$ is symmetric in $n$. Moreover, if $g$ is also symmetric with respect to permutations of its variables, then it suffices to obtain (2) for $\partial_\alpha= (\partial_1)^l.$

Let $$a_n:=|\widehat{g}(n)|\quad\mbox{ and}\quad\omega_N:=\omega(1/N),$$ where $\omega$ is a $r$-quasiconcave function with $\omega(t)\approx\omega_r(D^l g,t)_2$. Define $\mu$ as in Lemma \ref{lemma:discre}.
  
  For each $k\in I \cup J$, we construct a trigonometric polynomial $S_k$ with the following properties:
    \begin{enumerate}
        \item $\norm{S_k}_{C^l(\mathbb{T}^d)} \lesssim \omega_{\mu_k};$
        \item $\norm{S_k}_{C^{l+r}(\mathbb{T}^d)} \lesssim \omega_{\mu_{k+1}-1}(\mu_{k+1}-1)^r;$
        \item $\spec(S_k)\subset B_\infty(2(\mu_{k+1}-1)) \setminus B_\infty(\mu_k/2);$
        \item $|\widehat{S_k}(n)|\geq a_n, \, n \in B_\infty((\mu_{k+1}-1)) \setminus B_\infty(\mu_k-1);$
        \item $\widehat{S_k}(n)\in \mathbb{R}$.
       
    \end{enumerate}

    First, we assume that $k\in I$. Since \eqref{eq:L2disc} holds, by Theorem \ref{theorem:KKdL} with $s=l$ applied to the sequence $a$ supported on  $ B_\infty((\mu_{k+1}-1)) \setminus B_\infty(\mu_k-1)$, there exists a function $h \in C^l(\mathbb{T}^d)$ such that 
    
    \begin{itemize}
        \item $\widehat{h}(n)\in \mathbb{R};$
        \item   $|\widehat{h}(n)|\geq a_n ,  \qquad n \in B_\infty((\mu_{k+1}-1)) \setminus B_\infty(\mu_k-1);$
        \item $\norm{h}_{C^l(\mathbb{T}^d)} \lesssim \omega_{\mu_k}$.
       
    \end{itemize} Let $V^d_n$ be the $d$-dimensional de la Vallée-Poussin operator of degree $n$ defined in \eqref{eq:VP}. Observe that, for any $h\in L ^\infty$, one has \begin{itemize}
       \item[$(i)$]$\spec(V^d_n(h))\subset B_\infty(2n);$
     \item[$(ii)$] $\norm{V^d_n(h)}_{C^l(\mathbb{T}^d)} \lesssim \norm{h}_{C^l(\mathbb{T}^d)};$
   \item[$(iii)$] if $\spec(h)\subset[-n,n]^d$, then $V^d_n(h)=h$.
    \end{itemize}

   Let 
   \begin{align}\label{S_k}
       S_k:= V^d_{\mu_{k+1}-1}\left( h -V^d_{\mu_k/2}(h)\right).
        \end{align}The verification of items  (1), (3)--(5) is straightforward from the definition of $S_k$ and properties (i)--(iii) of $V^d_n$. For (2), note that by the Bernstein inequality (item (1) in Lemma \ref{lemma:bernstein}), \begin{align*}
       \norm{S_k}_{C^{l+r}(\mathbb{T}^d)} &\lesssim (\mu_{k+1}-1)^r \norm{S_k}_{C^l(\mathbb{T}^d)} \\&\lesssim  (\mu_{k+1}-1)^r\omega_{\mu_k} \approx (\mu_{k+1}-1)^r \omega_{\mu_{k+1}-1} ,
        \end{align*}
        where the last estimate follows from $k\in I$.

Second,    let $k\in J$. Taking into account \eqref{eq:L2disc2}, by Theorem \ref{theorem:KKdL} with $s=l+r$ applied to the sequence $a$ supported on $B_\infty(\mu_{k+1}-1) \setminus B_\infty(\mu_k-1)$, there exists a function $h$ such that 
    
    \begin{itemize}
        \item $\widehat{h}(n)\in \mathbb{R},$
        \item   $|\widehat{h}(n)|\geq a_n,  \qquad n \in B_\infty((\mu_{k+1}-1)) \setminus B_\infty(\mu_k-1);$
        \item $\norm{h}_{C^{l+r}(\mathbb{T}^d)} \lesssim \omega_{\mu_{k+1}-1} (\mu_{k+1}-1)^r$.
    \end{itemize} 
Define 
   $S_k$ as in 
   \eqref{S_k}.
   The verification of properties (2)--(5) is again clear. To see  (1), we note that by the reverse Bernstein inequality (item (2) in Lemma \ref{lemma:bernstein}), $$\mu^r_{k}\norm{S_k}_{C^l(\mathbb{T}^d)} \lesssim \norm{S_k}_{C^{l+r}(\mathbb{T}^d)} \lesssim \omega_{\mu_{k+1}-1} (\mu_{k+1}-1)^r \approx \mu^r_{k}\omega_{\mu_{k}} ,$$ where the last estimate follows from $k\in J$.

Finally, let
$$f=\sum_{k=1}^\infty  i ^ k S_k.$$ Note that, since $4\mu_k < \mu_{k+1}$, we have that $\spec(S_k) \cap \spec(S_{k+2})=\emptyset$. Thus, if $\norm{n}_\infty\in [\mu_k, \mu_{k+1}-1]$, then
$$|\widehat{f}(n)| \geq |\widehat{S_k}(n)+\widehat{S_{k-1}}(n)+\widehat{S_{k+1}}(n)| \geq  |\widehat{S_k}(n)| \geq a_n.$$ Note that this in particular implies that $$\omega_r(D^lf,N^{-1})_\infty \geq \omega_r(D^lf,N^{-1})_2 \gtrsim \omega_N.$$
 
After constructing $f$, we show that $\omega_r(D^lf,N^{-1})_\infty\lesssim \omega_N$. To this end,  let  $\mu_k\leq N <\mu_{k+1}$.
Then
\begin{align}
   \omega_r(D^lf,N^{-1})_\infty\leq  \left(\sum_{j<k} + \sum_{j=k} +\sum_{j>k}\right) \omega_r(D^lS_j,N^{-1})_\infty:=J_1+J_2+J_3.
\end{align}
   For $J_1$, recalling the properties of the $S_j$, we obtain the estimate \begin{equation}
       \label{eq:omeg1}\omega_r(D^lS_j,N^{-1})_\infty \lesssim \frac{1}{N^r}  \norm{S_j}_{C^{l+r}} \lesssim  \omega_{\mu_{j+1}-1} \frac{(\mu_{j+1}-1)^r}{N^r}, 
   \end{equation}  which, by the properties of $\mu$, implies
   \begin{align}
     J_1 \lesssim \frac{1}{N^r}  \sum_{j<k}  \omega_{\mu_{j+1}-1} (\mu_{j+1}-1)^r\lesssim \frac{1}{N^r} \omega_{\mu_k} \mu^r_k \leq \omega_N.
   \end{align}
    For $J_3$, using instead 
    \begin{equation}
        \label{eq:omeg2} \omega_r(D^lS_j,N^{-1})_\infty \lesssim  \norm{S_j}_{C^{l}} \lesssim  \omega_{\mu_{j}} 
    \end{equation}  we derive
   \begin{align}
    J_3 \lesssim  \sum_{j>k}  \omega_{\mu_{j}} \lesssim  \omega_{\mu_{k+1}} \leq \omega_N.
   \end{align}

   Finally, for $J_2$ there are two possibilities. First, if $k\in I$,
   we use \eqref{eq:omeg2} and obtain\begin{align}
       \omega_r(S_k^{(\alpha)},{N^{-1}} )_\infty\lesssim  \omega_{\mu_{k}} \approx \omega_N.
   \end{align}
  Second, if $k\in J$, we use \eqref{eq:omeg1} to obtain
   \begin{align}
        \omega_r(S_k^{(\alpha)},{N^{-1}} )_\infty\lesssim  \frac{1}{N^r}\omega_{\mu_{k+1}-1} (\mu_{k+1}-1)^r\approx \frac{\mu^r_{k}}{N^r}\omega_{\mu_{k}}  \approx\omega_N.
   \end{align}
   The proof is now complete.
\end{proof}

\begin{proof}[Proof of Theorem \ref{coro:hardy}]
 Theorem \ref{th:main} provides the equivalence between the  first 
 two items.
 
 Since Lemma \ref{lemma:l2} gives an expression for the norm of  a function in $
 {\operatorname{Lip}^{r,l}_2(\omega)}
 $ in terms of its Fourier coefficients, by using the change of variables $x_n=R_n n^l$ and observing that $$\sup_{a_k} \frac{\sum_{\norm{k}_\infty=n} |a_k|^p}{R_n^p} \approx n^{(d-1)(1-\frac p2)},$$ 
  we see that it suffices to characterize the $\omega_N$ for which there exists $K<\infty$ such that for any nonnegative sequence $(x_n)_{n=1}^\infty$
    \begin{equation}
           \label{ineq:xnK} \left(\sum_{n=1}^\infty x_n^p n^{(d-1)(1-\frac{p}{2})-pl} \right)^\frac{1}{p}\leq K\left(\sup_N \frac{1}{\omega^2_N} \sum_{n=1}^\infty \left(x_n \min(1, \frac{n}{N})^r\right)^2\right)^{\frac12}.\end{equation}

We now characterize \eqref{ineq:xnK}. To begin with, we treat some trivial cases. First, if $d(\frac1p-\frac12)\geq r+l$, we see that \eqref{ineq:xnK} cannot hold because it would imply the inequality
\begin{align*}
    \left(\sum_{n=1}^\infty x_n^p n^{(d-1)(1-\frac{p}{2})-pl} \right)^\frac{1}{p}
 & \leq K\left(\sup_N \frac{1}{\omega^2_N} \sum_{n=1}^\infty \left(x_n \min(1, \frac{n}{N})^r\right)^2\right)^{\frac12}
\\
 &\lesssim \left( \sum_{n=1}^\infty \left(x_nn^r\right)^2\right)^{\frac12},
 \end{align*}
which is not true for some $(x_n)$.

Second, if $d(\frac1p-\frac12)<l$, then
\begin{align*}\left(\sum_{n=1}^\infty x_n^p n^{(d-1)(1-\frac{p}{2})-pl} \right)^\frac{1}{p}
 &\lesssim \left(\sum_{n=1}^\infty x_n ^2\right)^\frac 12 \\
 &\lesssim \left(\sup_N \frac{1}{\omega^2_N} \sum_{n=1}^\infty \left(x_n \min(1, \frac{n}{N})^r\right)^2\right)^{\frac12},
 \end{align*}so \eqref{ineq:xnK} holds.

\iffalse
Third, if $-1< d-2l<2r$ and $\omega_N N^r $ is bounded, then
$$ \left(\sup_N \frac{1}{\omega^2_N} \sum_{n=1}^\infty \left(x_n \min(1, \frac{n}{N})^r\right)^2\right)^{\frac12} \approx \left( \sum_{n=1}^\infty \left(x_nn^r\right)^2\right)^{\frac12}.$$ In this case, we have that $K<\infty$ if and only if $\sum_{n=1}^\infty x_n  n^{\frac{d}{2}-l-1} \lesssim \left(\sum_{n=1}^\infty x_n^2 n^{2r}\right)^\frac 12$, which is implied by $d-2l<2r$.

\fi

Now assume that $r+l>d(\frac1p-\frac12)\geq l$. By \eqref{eq:L2disc} and \eqref{eq:L2disc2}, 
we know that

$$\left(\sup_N \frac{1}{\omega^2_N} \sum_{n=1}^\infty \left(x_n \min(1, \frac{n}{N})^r\right)^2\right)^{\frac12} \lesssim 1$$
holds for a nonnegative $(x_n)$
if and only if
   for all $k \in [1,L] $ we have
    \begin{align} 
    \label{eq:s1}
        \sum_{n=\mu_k}^{\mu_{k+1}-1} |x_n|^2 \lesssim& \omega^2_{\mu_k}\quad\mbox{ if }\quad k\in I \\
        \mbox{and}\\
        \label{eq:s2}
        \sum_{n=\mu_k}^{\mu_{k+1}-1} |n^{r}x_n|^2 \lesssim& (\mu_{k+1}-1)^{2r}\omega^2_{\mu_{k+1}-1} \quad \mbox{ if }\quad k\in J.
    \end{align}
Hence, 
 the best constant in \eqref{ineq:xnK} is given by
\begin{align}
K^p= & \sup_{x} \sum_{k\in I\cup J} \left(\sum_{n=\mu_k}^{\mu_{k+1}-1} x_n^p n^{(d-1)(1-\frac{p}{2})-pl}\right),
\end{align}
 where the supremum %in the first line 
  is taken over those $x$ which satisfy \eqref{eq:s1} and \eqref{eq:s2}. Then by 
Hölder's inequality,
\begin{align}
K^p&\approx  \sum_{k\in I} \omega^p_{\mu_k} \left(\sum_{n=\mu_k}^{\mu_{k+1}-1} n^{(d-1)- \frac{l}{\frac 1p-\frac12}}\right)^{1-\frac p2} \\
 &+ \sum_{k\in J} (\mu_{k+1}-1)^{pr}\omega^p_{\mu_{k+1}-1} \left(\sum_{n=\mu_k}^{\mu_{k+1}-1} n^{(d-1)-\frac{l+r}{\frac 1p-\frac12}}\right)^{1-\frac p2}\\
 &\approx  \sum_{k\in I \cup J} \left(\sum_{n=\mu_k}^{\mu_{k+1}-1} \omega^ \frac{1}{\frac 1p-\frac12}_n n^{d-1-\frac{l}{\frac 1p-\frac12}}\right)^{1-\frac p2},
\end{align}
where the last equivalence follows from Lemma \ref{lemma:discre} (6).

An application of Lemma \ref{theorem:split} with $$q =1- \frac{p}{2},\; \; \; \overline \omega^p_N=  \sum_{n=1}^\infty \frac{\alpha_n}{n^{pr}+N^{pr}},\quad \mbox{and }\quad f_n = n^{d-1-\frac{l}{\frac 1p-\frac12}},$$ where $\overline \omega$ is obtained from $\omega$ by using Lemma \ref{lemma:majorant}, yields
$$K^p\approx \sum_{n=1}^ \infty \alpha_n \left(\sum_{j=1}^\infty \frac{j^{d-1-\frac{l}{\frac 1p-\frac12}}}{j^{\frac{r}{\frac 1p-\frac12}} +n^{\frac{r}{\frac 1p-\frac12}}}\right)^{1-\frac p2}.$$
Observe that in the range $r+l>d(\frac1p-\frac12)\geq l$,
$$\sum_{j=1}^\infty \frac{j^{d-1-\frac{l}{\frac 1p-\frac12}}}{j^{\frac{r}{\frac 1p-\frac12}} +n^{\frac{r}{\frac 1p-\frac12}}}\approx\begin{cases}
    n^{d-\frac{l+r}{\frac 1p-\frac12}}, \quad d> \frac{l}{\frac 1p-\frac12},\\
 n^{-\frac{r}{\frac 1p-\frac12}}     \log (n+1)\,,\quad  d=\frac{l}{\frac 1p-\frac12}.
\end{cases}$$

In the former case, by repeated integration by parts we obtain
\begin{align*}
K^p \approx \sum_{n=1}^ \infty \alpha_n n^{p ( d( \frac1p- \frac12)-l-r)} &\approx \sum_{n=1}^\infty n^{p(d(\frac1p -\frac12)-l)-1}  \left(\sum_{m=n}^\infty m^{-pr-1} \sum_{l=1}^m \alpha _l \right) \\& \approx \sum_{n=1}^\infty n^{p(d(\frac1p -\frac12)-l)-1}\omega^p_n,
\end{align*}
where the last estimate follows from
$$\sum_{m=n}^\infty m^{-pr-1} \sum_{l=1}^m \alpha _l \approx \sum_{m=1}^\infty \frac{\alpha_m}{n^{pr}+m^{pr}}\approx \omega^p_n.$$

Thus, in this case 
$$K^p\approx \int_0 ^1 \omega^p(t) t^{p(-d(\frac1p -\frac12)+l)} \frac{dt}{t}.$$
In the latter case, 
\begin{align*}
    K^p\approx \sum_{n=1}^ \infty \alpha_n \frac{\log^{1 -\frac p2} (n+1)}{n^{pr}}
&\approx  \sum_{n=1}^\infty \frac {1}{n \log^\frac p2 (n+1) }\left(\sum_{m=n}^\infty m^{-pr-1} \sum_{l=1}^m \alpha _l \right)  
\\&\approx \sum_{n=1}^\infty \frac{1}{n \log^\frac p2 (n+1) } \omega^p_n. 
\end{align*}
In other words, we derive that
$$K^p \approx \int_0^1 \omega^p(t) \Big(\log \frac2t\Big)^{-\frac p2} \frac{dt}{t},$$
completing the  proof.
\end{proof}

\medskip
\section{Absolute integrability of Fourier transforms}

%Finally, for the Fourier transform
%$$ \widehat{f}(\xi)=\int_{\mathbb{R}^d} {f}(x) e^{-2\pi i \langle x, \xi \rangle} d x,
%$$we obtain an identical relation between $L^2$ smoothness and absolute integrability of $\widehat{f}$.

Historically, 
the study of 
absolute integrability of the Fourier transform
$$ \widehat{f}(\xi)=\int_{\mathbb{R}^d} {f}(x) e^{-2\pi i \langle x, \xi \rangle} d x
$$
was initiated   in 1927
by
Titchmarsh \cite{tit}, who found a Bernstein-type condition given in terms of the decay of the $L_q$-modulus of smoothness of  $f$. 
 In other words, this condition ensures that a function belongs to the space
 $$\mathcal{A}_p(\mathbb{R}^d) :=
\Big\{ f\in L ^2(\mathbb{R}^d): \widehat{f}\in L^p (\mathbb{R}^d)
\Big\}.$$
A slightly different definition of 
the 
Wiener space was suggested by
Beurling  \cite{Beurling}. %, where he studied % obtained first results on

There are many known sufficient and some necessary conditions for smooth classes to belong to $\mathcal{A}_p(\mathbb{R}^d)$. For further details, we refer to the survey \cite{samko} and the monograph \cite[Chapter 6]{peetre}.
Similarly to the periodic case, we obtain necessary and sufficient conditions in terms of the $L_2$-moduli of smoothness.
(The space ${\operatorname{Lip}^{r,l}_2(\omega;\mathbb{R}^d)}$ 
is defined as ${\operatorname{Lip}^{r,l}_2(\omega;\mathbb{T}^d)}$ in the introduction 
replacing 
 $\sup _{0<t<1} $ by $\sup _{0<t<\infty}$.)
\begin{theorem} Let $l\in \mathbb{N}_0$,
$r\in \mathbb{N}$, and $\omega$ be $r$-quasiconcave. Let $0<p<2$. Then, the following are equivalent:
    \label{prop:hardy}
\begin{enumerate}
      \item $$L^2(\mathbb{R}^d) \cap
      {\operatorname{Lip}^{r,l}_2(\omega;\mathbb{R}^d)}
            \subset \mathcal{A}_p(\mathbb{R}^d);$$
   \item one of the following holds:
    \begin{itemize}
       \item[$(i)$] $d(\frac 1p- \frac12)<l;$
       \item[$(ii)$] $l<d(\frac 1p- \frac12)<l+r$ and $$\int_0 ^1 \omega^p(t) t^{p(-d(\frac1p -\frac12)+l)} \frac{dt}{t}<\infty;$$
        \item[$(iii)$] $d(\frac 1p- \frac12)=l$ and $$\int_0 ^1 \frac{\omega^p(t)}{\big(\log \frac2t\big)^{\frac{p}{2}}} \frac{dt}{t}<\infty.$$
        \end{itemize}
\end{enumerate}
\end{theorem}

\begin{proof}%[Proof of Proposition \ref{prop:hardy}]

Analogously to the discrete case, we need to characterize the majorant $\omega$ for which 

\begin{equation}
\label{eq:step0}
    \left(\int_0^\infty g^2(t) dt  \right)^\frac12+\left(\sup _{0<h<\infty} \frac{1}{\omega^2(h)} \int_0 ^\infty \left(g(t)t^l \min(1,th)^r\right)^2 dt\right)^{\frac12}<\infty
\end{equation} implies

\begin{equation}
\label{eq:step}
    \int_0 ^\infty g^p(t) t^{(d-1)(1-\frac{p}{2})} dt<\infty.
\end{equation}
    Since \eqref{eq:step0} implies that $g\in L^2$, 
    \eqref{eq:step} holds if and only if
    \begin{equation}
\label{eq:step2}
    \int_1 ^\infty g^p(t) t^{(d-1)(1-\frac{p}{2})} dt<\infty.
\end{equation}
Thus, instead of \eqref{eq:step0}, we may assume that
 \begin{align}
     I:=\left(\int_1^\infty g^2(t) dt  \right)^\frac12+\left(\sup _{0<h<\infty} \frac{1}{\omega^2(h)} \int_1 ^\infty \left(g(t)t^l \min(1,th)^r\right)^2 dt\right)^{\frac12}<\infty.
 \end{align}
 Since \begin{equation}
     I\approx  \left(\sup _{0<h<1} \frac{1}{\omega^2(h)} \int_1 ^\infty \left(g(t)t^l \min(1,th)^r\right)^2 dt\right)^{\frac12},\end{equation}
 it suffices to study (cf. \eqref{ineq:xnK}) the following inequality for non-negative $g$:
    \begin{align}
       &\left( \int_1^\infty g^p(t) t^{(d-1)(1-\frac{p}{2})-pl} dt\right)^\frac 1p\\&\leq K\left(\sup _{0<h<1} \frac{1}{\omega^2(h)} \int_1 ^\infty \left(g(t) \min(1,th)^r\right)^2 dt\right)^{\frac12},
    \end{align} which is characterized 
    just like inequality \eqref{ineq:xnK}.
\end{proof}

\medskip

\bibliographystyle{amsplain}
\bibliography{main}

\end{document}